\documentclass{emsprocart}

\usepackage{mathrsfs,amssymb,amsfonts,amsmath}

\usepackage{graphics}
\usepackage[all]{xy}
\xyoption{curve}
\xyoption{import}
\xyoption{arc}
\xyoption{ps}
\UsePSspecials{dvips}
\usepackage{amsmath}
\usepackage{graphicx}
\usepackage[all]{xy}

\contact[carrillo@math.jussieu.fr]{Paulo Carrillo Rouse\\ Projet
  d'alg{\`e}bres d'op{\'e}rateurs\\ Universit{\'e} de Paris 7\\ 175, rue de
  Chevaleret\\ 
Paris, France}

\newtheorem{theorem}{Theorem}[section]

\newtheorem{proposition}[theorem]{Proposition}

\theoremstyle{definition}

\newtheorem{definition}[theorem]{Definition}

\newtheorem{remark}[theorem]{Remark}

\newcommand{\ci}{C^{\infty}}
\newcommand{\Cat}{\mathscr{C}}
\newcommand{\Dnc}{\mathscr{D}}
\newcommand{\gr}{\mathscr{G}}
\newcommand{\go}{\mathscr{G} ^{(0)}}

\newcommand{\Nb}{\mathscr{N}}
\newcommand{\sw}{\mathscr{S}}
\newcommand{\Uo}{\mathscr{U}}
\newcommand{\Vo}{\mathscr{V}}
\newcommand{\Rr}{\mathbf{R}}
\newcommand{\Nat}{\mathbf{N}}
\newcommand{\src}{\mathscr{S}_{r,c}}
\newcommand{\cg}{C_{c}^{\infty}(\gr)}
\newcommand{\cgo}{C_{c}^{\infty}(\go)}

\newcommand{\ckt}{C_{c}^{k}(\gr \times [0,1])}
\newcommand{\ck}{C_{c}^{k}(\gr)}

\title[An analytic index for Lie groupoids]{An analytic 
index for Lie groupoids}

\author{Paulo Carrillo Rouse}

\begin{document}

\begin{abstract}
For a Lie groupoid there is an analytic index morphism 
which takes values in the
$K-$theory of the $C^*$-algebra associated to the groupoid. This 
is a good invariant but extracting numerical invariants from it, 
with the existent
tools, is very difficult. 
In this work, we define another analytic index morphism associated
to a Lie groupoid; this one takes
values in a group that allows us to do pairings with cyclic
cocycles. This last group is related to the
compactly supported functions on the groupoid. We use the tangent
groupoid to define our index as
a sort of ''deformation''.
\end{abstract}

\begin{classification}

Primary 19-06; Secondary 19K56.

\end{classification}

\begin{keywords}

Lie groupoids, Tangent groupoid, K-theory, Index theory.

\end{keywords}

\maketitle

\section{Introduction}
The concept of groupoid is central in non commutative geometry. Groupoids
generalize the concepts of spaces, groups and 
equivalence relations. In the late 70's, mainly with the work
of Alain Connes, it became clear that groupoids are natural
substitutes of singular spaces. Furthermore, Connes showed that many
groupoids and algebras associated to them appeared as `non commutative
analogues` of smooth manifolds to which many tools of geometry such as
K-theory and Characteristic classes could be applied.\\
One classical way to obtain invariants is through the index
theory in the sense of Atiyah-Singer. 
Given $\gr \rightrightarrows \go$ a Lie groupoid, it is possible to
talk about Pseudodifferential calculus on it
(See \cite{Co1},\cite{CS},\cite{MP}, \cite{NWX}). Hence, if one
has a $\gr$-pseudodifferential operator $P$, then one can consider two
elements associated to it
\begin{itemize}
\item a symbol  class $[\sigma_P]\in K^0(A^*\gr)$,
\item an index $ind (P)\in K_0(\ci_c (\gr))$.
\end{itemize}   
In fact, the index of a $\gr$-pseudodifferential operator 
is usually considered as an element of
$K_0(C^*(\gr))$ (in the present work
$C^*(\gr)$ stands for the reduced $C^*$-algebra) rather 
than an element of $K_0(\ci_c (\gr))$. The idea is to 
take the image of $ind(P)$ through the $K$-theory morphism
\begin{align}\label{j} 
K_0(\ci_c (\gr)) \stackrel{j}{\rightarrow} K_0(C^*(\gr))  
\end{align}
induced by the inclusion $\ci_c (\gr) \hookrightarrow C^*(\gr)$. Such an 
image, noted by $ind_a(P)$, is called "The analytic
index of the Lie groupoid $\gr$". The reason for considering the
"$C^*$-index" $ind_a(P)$, is that it is a homotopy invariant of the  
$\gr$-pseudodifferential elliptic operators while the index $ind(P)$
is not. Even more, the analytic index can be
given as a group morphism (using the Connes tangent groupoid for example)
$$ 
ind_a:K^0(A^*\gr) \rightarrow K_0(C^*(\gr)) \nonumber
$$
that fits in the following commutative diagram:
\[
\xymatrix{
\{\gr-pdo ell.\} \ar[dd]_{symbol} \ar[rr]^{ind} & &
K_0(\ci_c(\gr)) \ar[dd]^j\\ & & \\
K_0(A^*\gr) \ar[rr]_{ind_a} &  & K_0(C_{r}^{*}(\gr))
}
\]
We cannot expect (\ref{j}) to be
injective even if it comes from the canonical inclusion. 
In \cite{Co2}, Alain Connes gives examples of the non injectivity of (\ref{j}) and he
discusses other reasons why it is not sufficient, in general, to stay with the
index at the $\ci_c$-level. The analytic index $ind_a$ 
is a good invariant, however, extracting numerical invariants from it, with the existent tools, is very difficult. Moreover, roughly speaking, the fact that
the morphism in (\ref{j}) is not injective tells us
that in some way we are loosing information.\\ 

In the present work we propose
an intermediate group $K_0(\ci_c(\gr)) \rightarrow
 K_{0}^{h,k}(\gr)\rightarrow K_0(C^*(\gr))$ (one for each
$k \in \Nat$) and we construct "analytic index morphism" 
$$
ind_{a}^{h,k}:K_0(A^*\gr) 
\rightarrow K_{0}^{h,k}(\gr)$$ 
that will allow us to obtain more explicit invariants. Let us briefly explain 
what is this group $K_{0}^{h,k}(\gr)$, why we call $ind_{a}^{h,k}$ an "analytic index", and
why we state that it will allow to obtain 
more explicit invariants.
\begin{itemize}
\item[(a)]The group $K_{0}^{h,k}(\gr)$ is simply
the group $K_0(C_{c}^{k}(\gr))$ modulo the following homotopy equivalence
relation:\\
Given $x,y\in K_0(C_{c}^{k}(\gr))$ we say that $x$ and $y$ are homotopic, 
$x\sim_hy$, if and only if there exists $z\in K_0(C_{c}^{k}(\gr \times
[0,1]))$ such that $e_0(z)=x$ and $e_1(z)=y$ where $e_t$ is the
morphisms in $K$-theory induced by the respectives evaluations. One proves that this gives an equivalence relation and 
$K_0(C_{c}^{k}(\gr))/\sim_h$ is a group under the natural operation 
$[x]_h + [y]_h=[x+y]_h$.
\item[(b)]We call $ind_{a}^{h,k}$ an "analytic index " because it fits in the following
 commutative
diagram:

\[
\xymatrix{
\{\gr-pdo ell.\} \ar[dd]_{symb} \ar[r]^{ind}&
K_0(\ci_c(\gr)) \ar[d] & \\ & K_{0}^{h,k}(\gr) \ar[d] & \\
K_0(A^*\gr) \ar[r]_{ind_a} \ar[ru]^{ind_{a}^{h,k}}& K_0(C_{r}^{*}(\gr)) & 
}
\]
\item[(c)]We would like to obtain information from
  $K-$theory elements through the "Chern-Weil theories". In the
  classical case of spaces (manifolds, locally compact spaces, etc.) this is done by a Chern character. Now, Connes
  developped in the early 80's a way to generalize this to more general
  "spaces". In particular he introduced the Cyclic cohomology for general
  algebras and he defined 
pairings beetween $K-$theory elements and cyclic cocycles. There is
  also the notion of Chern character for the non
  commutative case.\\
Now, the cyclic homology theory is only
  interesting if it is applied to "smooth" algebras, typically
  algebras as $\ci(M)$ (for $M$ a manifold) or $\ci_c(\gr)$. In our case, it is possible
  to pair cyclic cocycles with the index $ind(P)\in
  K_0(\ci_c(\gr))$ but it is not always possible to extend this action
  to the $K_0(C^*(\gr))$.\\
What we are proposing has the advantages of both indices ($ind$
  and $ind_a$), in one hand it is a well defined group morphism and it
  defines a homotopy invariant, and on the other hand it is
  possible to extract numerical invariants using the existent tools. 
\end{itemize}
The article is organized as follows:\\
In the second section we recall the basic facts about Lie groupoids. We
explain very briefly how to define the index $ind(P) \in
K_0(\cg)$. For a complete exposition see \cite{NWX}. 
In the third section we explain the ``deformation to the normal cone``
construction associated to an injective immersion. Even if this could
be considered as classical material we do it in some detail since we
will use in the sequel very explicit descriptions. A particular case
of this construction is the tangent groupoid associated to a Lie groupoid. This
last example will be fundamental in our construction of the index.\\
In the fourth section we will construct an algebra of functions over
the tangent groupoid that satisfy a rapid decay condition at zero and
are compactly supported elsewhere. This algebra is the main point of this
work and its $K-$theory allows us to define the index
$ind_{a}^{h,k}$ as a "deformation" (in
the same spirit as in \cite{Co2},\cite{HS},\cite{MP}). 
The last section is devoted to the construction of the analytic
index.\\
All the results of the present work are part of the author's PHD thesis. 
Complete proofs and details may be found in \cite{Ca} and will be
published elsewhere.\\ 

\thanks{I want to thank my PHD advisor, Georges Skandalis, for
  all the ideas that he shared with me. I would like also to thank him for all the comments and remarks he made to the present work.}

\section{Lie groupoids}

Let us recall what a groupoid is:

\begin{definition}
A $\it{groupoid}$ consists of the following data:
two sets $\gr$ and $\go$, and maps
\begin{itemize}
\item[$\cdot$]$s,r:\gr \rightarrow \go$ 
called the source and target map respectively,
\item[$\cdot$]$m:\gr^{(2)}\rightarrow \gr$ called the product map\\ 
(where $\gr^{(2)}=\{ (\gamma,\eta)\in \gr \times \gr : s(\gamma)=r(\eta)\}$),
\item[$\cdot$]$u:\go \rightarrow \gr$ the unit map and 
\item[$\cdot$]$i:\gr \rightarrow \gr$
  the inverse map
\end{itemize}
such that, if we note $m(\gamma,\eta)=\gamma \cdot \eta$, $u(x)=x$ and 
$i(\gamma)=\gamma^{-1}$, we have 
\begin{itemize}
\item[1.]$\gamma \cdot (\eta \cdot \delta)=(\gamma \cdot \eta )\cdot \delta$, 
$\forall \gamma,\eta,\delta \in \gr$ when this is possible.
\item[2.]$\gamma \cdot x = \gamma$ and $x\cdot \eta =\eta$, $\forall
  \gamma,\eta \in \gr$ with $s(\gamma)=x$ and $r(\eta)=x$.
\item[3.]$\gamma \cdot \gamma^{-1} =u(r(\gamma))$ and 
$\gamma^{-1} \cdot \gamma =u(s(\gamma))$, $\forall \gamma \in \gr$.
\item[4.]$r(\gamma \cdot \eta) =r(\gamma)$ and $s(\gamma \cdot \eta) =s(\eta)$.
\end{itemize}

Generally, we denote a groupoid by $\gr \rightrightarrows \go $ where 
the parallel arrows are the source and target maps
and the other maps are given.

\end{definition}

Now, a Lie groupoid is a groupoid in which every set and map involved
in the last definition is $\ci$ (possibly with borders), and the source and target maps are
submersions. For $A,B$ subsets of $\go$ we use the notation
$\gr_{A}^{B}$ for the subset $\{ \gamma \in \gr : s(\gamma) \in A,\, 
r(\gamma)\in B\}$.\\
All along this paper, $\gr \rightrightarrows \go $ is going to be a
Lie groupoid.\\ 
We recall how to define an algebra structure in $\cg$ using
smooth Haar systems.
 
\begin{definition}
A $\it{smooth\, Haar\, system}$ over a Lie groupoid consists of a family of
measures $\mu_x$ in $\gr_x$ for each $x\in \go$ such that,

\begin{itemize}
\item for $\eta \in \gr_{x}^{y}$ we have the following compatibility
  condition:
$$\int_{\gr_x}f(\gamma)d\mu_x(\gamma)
=\int_{\gr_y}f(\gamma \circ \eta)d\mu_y(\gamma)$$
\item for each $f\in \cg$ the map

$$x\mapsto \int_{\gr_x}f(\gamma)d\mu_x(\gamma) $$ belongs to $\cgo$

\end{itemize}

\end{definition}

A Lie groupoid always has a smooth Haar system. In fact, if we
fix a smooth (positive) section of the 1-density bundle associated to
the Lie algebroid we obtain a smooth Haar system 
in a canonical way. The advantage of using 1-densities is
that the measures are locally equivalent to the Lebesgue measure. 
We suppose for the rest of tha paper that the smooth 
Haar systems are given by 1-densities (for complete
details see \cite{Pat}). 
We can now define a convolution
product on $\cg$: Let $f,g\in \cg$, we set

$$(f*g)(\gamma)
=\int_{\gr_{s(\gamma)}}
f(\gamma \cdot \eta^{-1})g(\eta)d\mu_{s(\gamma)}(\eta)$$

This gives a well defined associative product. 
\begin{remark}
There
is a way to avoid Haar systems when one works with Lie groupoids,
using the half densities (see Connes book \cite{Co2}).
\end{remark}
As we have 
mentioned in the introduction, we are going to consider some elements
in the $K$-theory group $K_0(\cg)$. We recall how this elements are
usually defined (See \cite{NWX} for complete details):\\
First we recall what a $\gr$-Pseudodifferentiel operator is:
\begin{definition}[$\gr$-PDO]
A $\gr$-$\it{Pseudodifferential}$ $\it{operator}$ is a familiy of
pseudodifferential operators $\{ P_x\}_{x\in \go} $ acting in
$\ci_c(\gr_x)$ such that if $\gamma \in \gr $ and
$$U_{\gamma}:\ci_c(\gr_{s(\gamma)}) \rightarrow \ci_c(\gr_{r(\gamma)}) $$
the induced operator, then we have the following compatibility condition 
$$ P_{r(\gamma)} \circ U_{\gamma}= U_{\gamma} \circ P_{s(\gamma)}$$
There is also a differentiability condition with respect to $x$ that can
be found in \cite{NWX}.

\end{definition}

Let $P$ be a $\gr$-Pseudodifferential elliptic operator. By definition
this means that 
there exists a parametrix, $\it{i.e.}$, a $\gr$-Pseudodifferential
 operator $Q$ such that $PQ-1$ and $QP-1$ belong to $\cg$
(where we are identifing $\cg$ with $\Psi^{-\infty}(\gr)$ as in
\cite{NWX}). The last data defines a quasi-isomorphism in 
$(\Psi^{+\infty},\cg)$ and then an element in $K_0(\cg)$ that we call
the index $ind(P)$. Similarly to the classical case, a $\gr$-PDO
operator has a principal symbol that defines an element in the
$K$-theory group $K^0(A^*\gr)$. 

\section{Deformation to the normal cone}

Let $M$ be a $\ci$ manifold and $X\subset M$ be a $\ci$ submanifold. We denote
by $\Nb_{X}^{M}$ the normal bundle to $X$ in $M$, $\it{i.e.}$, 
$\Nb_{X}^{M}:= T_XM/TX$.\\
We define the following set
\begin{align}
\Dnc_{X}^{M}:= \Nb_{X}^{M} \times {0} \bigsqcup M \times (0,1] 
\end{align} 
The purpose of this section is to recall how to define a $\ci$-structure with
boundary in $\Dnc_{X}^{M}$. This is more or less clasical, for example
it was extensively used in \cite{HS}. Here we are only going to
do a sketch.\\
Let us first consider the case where $M=\Rr^n$ 
and $X=\Rr^p \times \{ 0\}$ (where we
identify canonically $X=\Rr^p$). We denote by
$q=n-p$ and by $\Dnc_{p}^{n}$ for $\Dnc_{\Rr^p}^{\Rr^n}$ as above. In this case
we clearly have that $\Dnc_{p}^{n}=\Rr^p \times \Rr^q \times [0,1]$ (as a
set). Consider the 
bijection  $\Psi: \Rr^p \times \Rr^q \times [0,1] \rightarrow
\Dnc_{p}^{n}$ given by 
$$
\Psi(x,\xi ,t) = \left\{ 
\begin{array}{cc}
(x,\xi ,0) &\mbox{ if } t=0 \\

(x,t\xi ,t) &\mbox{ if } t>0
\end{array}\right.
$$
which inverse is given explicity by 
$$
\Psi^{-1}(x,\xi ,t) = \left\{ 
\begin{array}{cc}
(x,\xi ,0) &\mbox{ if } t=0 \\
(x,\frac{1}{t}\xi ,t) &\mbox{ if } t>0
\end{array}\right.
$$
We can consider the $\ci$-structure with border on $\Dnc_{p}^{n}$
induced by this bijection.\\  

In the general case, let 
$(\Uo,\phi)$ be a local chart in $M$ and suppose it is an $X$-slice,
so that it
satisfies

\begin{itemize}
\item[1)]$\phi : \Uo \stackrel{\cong}{\rightarrow} U \subset \Rr^p\times \Rr^q$
\item[2)]If $\Uo \cap X =\Vo$, $\Vo=\phi^{-1}( U \cap \Rr^p \times \{ 0\}
  )$ (we note $V=U \cap \Rr^p \times \{ 0\}$)
\end{itemize}
With this notation we have that $\Dnc_{V}^{U}\subset \Dnc_{p}^{n}$ is an
open subset. We may define a function $$
\tilde{\phi}:\Dnc_{\Vo}^{\Uo} \rightarrow \Dnc_{V}^{U} 
$$ in the following way: For $x\in \Vo$ we have $\phi (x)\in \Rr^p
\times \{0\}$. If we write 
$\phi(x)=(\phi_1(x),0)$, then 
$$ \phi_1 :\Vo \rightarrow V \subset \Rr^p$$ 
is a diffeomorphism, where $V=U\cap (\Rr^p \times \{0\})$. We set 
$\tilde{\phi}(v,\xi ,0)= (\phi_1 (v),d_N\phi_v (\xi ),0)$ and 
$\tilde{\phi}(u,t)= (\phi (u),t)$ 
for $t\neq 0$. Here 
$d_N\phi_v: N_v \rightarrow \Rr^q$ is the normal component of the
 derivate $d\phi_v$ for $v\in \Vo$. It is clear that $\tilde{\phi}$ is
 also a  bijection (in particular it induces a $C^{\infty}$ structure 
with border over $\Dnc_{\Vo}^{\Uo}$).\\ 
Now, let us consider an atlas 
$ \{ (\Uo_{\alpha},\phi_{\alpha}) \}_{\alpha \in \Delta}$ of $M$
 consisting of $X-$slices. Then we have the following proposition: 

\begin{proposition}\label{atlas}
The collection $ \{ (\Dnc_{\Vo_{\alpha}}^{\Uo_{\alpha}},\tilde{\phi_{\alpha})}
  \} _{\alpha \in \Delta }$ is a $\ci$-atlas with border over
  $\Dnc_{X}^{M}$.
\end{proposition}

\begin{definition}[DNC]
Let $X\subset M$ be as above. The set
$\Dnc_{X}^{M}$ provided with the  $C^{\infty}$ structure with border
induced by the atlas described in the last proposition is called
$\it{``The\, deformation\, to\, the\, normal\, cone\, associated\, to\,}$   
$X\subset M$``. We will often write DNC instead of
Deformation to the normal cone. 
\end{definition}

\begin{remark}
Following the same steps, we can define the deformation to the normal
cone associated to an injective immersion $X\hookrightarrow M$.
\end{remark}

The most important feature about the DNC construction is that it is in
some sense functorial. More explicity, let $(M,X)$ 
and $(M',X')$ be $\ci$-couples as above and let
 $F:(M,X)\rightarrow (M',X')$
be a couple morphism, i.e., a $\ci$ map   
$F:M\rightarrow M'$, with $F(X)\subset X'$. We define 
$ \Dnc(F): \Dnc_{X}^{M} \rightarrow \Dnc_{X'}^{M'} $ by the following formulas:\\

$\Dnc(F) (x,\xi ,0)= (F(x),d_NF_x (\xi),0)$ and\\

$\Dnc(F) (m ,t)= (F(m),t)$ for $t\neq 0$,\\

where $d_NF_x$ is by definition the map
\[  (\Nb_{X}^{M})_x 
\stackrel{d_NF_x}{\longrightarrow}  (\Nb_{X'}^{M'})_{F(x)} \]
induced by $ T_xM 
\stackrel{dF_x}{\longrightarrow}  T_{F(x)}M'$.\\
We have the following proposition.
\begin{proposition}
The map $\Dnc(F):\Dnc_{X}^{M} \rightarrow \Dnc_{X'}^{M'}$ is $\ci$.
\end{proposition}

\begin{remark}
If we consider the category $\Cat_{2}^{\infty}$ of  $\ci$ pairs given by
a $\ci$ manifold and a $\ci$ submanifold, and pair morphisms as above,
we can reformulate the proposition and say that we have a functor
$$\Dnc : \Cat_{2}^{\infty} \rightarrow \Cat^{\infty}$$ where
$\Cat^{\infty}$ denote the category of  $\ci$ manifolds with border.\\
\end{remark}

\subsection{The tangent groupoid}

\begin{definition}[Tangent groupoid]
Let $\gr \rightrightarrows \go $ be a Lie groupoid. $\it{The\, tangent\,
groupoid}$ associated to $\gr$ is the groupoid that has $\Dnc_{\go}^{\gr}
$ as the set of arrows and  $\go \times [0,1]$ as the units, with:\\

\begin{itemize}

\item[$\cdot$] $s^T(x,\eta ,0) =(x,0)$ and $r^T(x,\eta ,0) =(x,0)$ at $t=0$.
\item[$\cdot$] $s^T(\gamma,t) =(s(\gamma),t)$ and $r^T(\gamma,t)
  =(r(\gamma),t)$ at $t\neq0$.
\item[$\cdot$] The product is given by
  $m^T((x,\eta,0),(x,\xi,0))=(x,\eta +\xi ,0)$ et \linebreak $m^T((\gamma,t), 
  (\beta ,t))= (m(\gamma,\beta) , t)$ if $t\neq 0 $ and 
if $r(\beta)=s(\gamma)$.
\item[$\cdot$] The unit map $u^T:\go \rightarrow \gr^T$ is given by
 $u^T(x,0)=(x,0)$ and $u^T(x,t)=(u(x),t)$ for $t\neq 0$.
\end{itemize}
We denote $\gr^{T}:= \Dnc_{\go}^{\gr}$.

\end{definition} 

As we have seen above, $\gr^{T}$ can be considered as a $\ci$ manifold with
border. As a consequence of the functoriality of the DNC construction
we can show that the tangent groupoid is in fact a Lie
groupoid. Indeed, it is easy to check that if we identify in a
canonical way $\Dnc_{\go}^{\gr^{(2)}}$ with $(\gr^T)^{(2)}$, then 
$$ m^T=\Dnc(m),\, s^T=\Dnc(s), \,  r^T=\Dnc(r),\,  u^T=\Dnc(u)$$
where we are considering the following pair morphisms:
\begin{align}  
m:((\gr)^{(2)},\go)\rightarrow (\gr,\go ), \nonumber
\\
s,r:(\gr ,\go) \rightarrow (\go,\go),\nonumber 
\\
u:(\go,\go)\rightarrow (\gr,\go ).\nonumber
\end{align}
Finally, if $\{ \mu_x\}$ is a smooth Haar system on $\gr$, then, setting

\begin{itemize}
\item $\mu_{(x,0)}:=\mu_x$ at $(\gr^T)_{(x,0)}=T_x\gr_x$ and
\item $\mu_{(x,t):=t^{-q}\cdot \mu_x} $ at $(\gr^T)_{(x,t)}=\gr_x$ for
  $t\neq 0$, where $q=dim\, \gr_x$, 
\end{itemize}
one obtains a smooth Haar system for the Tangent groupoid 
(details may be found in \cite{Pat}).
\section{An algebra for the Tangent groupoid}

In this section we will show how to construct an algebra for the
tangent groupoid which consist of $\ci$ functions that satisfy a rapid
decay condition at zero while out of zero they satisfy a compact
support condition. This algebra is the main construction in
this work.

\subsection{Schwartz type spaces for Deformation to the normal cone manifolds}

Our algebra for the Tangent groupoid will be a particular
case of a construction associated to any deformation to the normal
cone.\\
Before giving the general definition, let us start locally. Let $p, q\in \Nat$ 
and $U \subset \Rr^p
  \times \Rr^q$ 
an open subset, and let $V=U\cap (\Rr^p \times \{ 0\})$. We set 
$$ \Omega_{V}^{U}:=\{ (x,\xi,t)\in \Rr^p
  \times \Rr^q \times [0,1] : (x,t\cdot \xi)\in U\},$$
which is an open subset of $\Rr^p
  \times \Rr^q \times [0,1]$ and is diffeomorphic to $\Dnc_{V}^{U}$ via the restriction of the map $\Psi$ used to define the $\ci$ structure on $\Dnc_{p}^{n}$ (where, as above, $n=p+q$). We can now give the following definition.
\begin{definition}\label{ladef}
Let $p, q\in \Nat$ and $U \subset \Rr^p
  \times \Rr^q$ 
an open subset, and let $V=U\cap (\Rr^p \times \{ 0\})$.
\begin{itemize}
\item[(1)]Let $K\subset U \times [0,1]$ be a compact
  subset. We say that $K$ is a $\it{conic\, compact}$ subset of $U \times [0,1]$
relative to $V$ if
\[ K_0=K\cap (U \times \{ 0\}) \subset V\]

\item[(2)]Let $g \in \ci
  (\Omega_{V}^{U})$. We say that
   $f$ has $\it{conic\, compact\, support}$ $K$, if there exists a conic
  compact $K$
 of $U \times [0,1]$ relative to $V$ such that if $t\neq 0$ and 
$(x, t\xi ,t) \notin K$ then $g(x, \xi ,t)=0$.

\item[(3)]We denote by $\src (\Omega_{V}^{U})$ 
the set of functions
$g\in \ci (\Omega_{V}^{U})$ 
that have compact conic support and that satisfy the following condition:

\begin{itemize}
\item[$(s_1$)]$\forall$ $k,m\in \Nat$, $l\in \Nat^p$
and $\alpha \in \Nat^q$ there exists $C_{(k,m,l,\alpha)} >0$ such that

\[ (1+\| \xi \|^2)^k \| \partial_{x}^{l}\partial_{\xi}^{\alpha}
\partial_{t}^{m}g(x,\xi ,t) \| \leq C_{(k,m,l,\alpha)}   \]
\end{itemize}

\end{itemize}

\end{definition}

Now, the spaces $\src (\Omega_{V}^{U})$ are invariant under
diffeomorphisms. More precisely, let $F:U\rightarrow U'$ be a $\ci$ diffeomorphism where 
$U\subset \Rr^p \times \Rr^q$ and $U'\subset \Rr^{p} \times \Rr^{q}$ are
open subsets. We write $F=(F_1,F_2)$ and we suppose that
$F_2(x,0)=0$. Then the function 
$\tilde{F}:\Omega_{V}^{U} \rightarrow \Omega_{V'}^{U'}$ defined by
$$
\tilde{F}(x,\xi ,t) = \left\{ 
\begin{array}{cc}
(F_1(x,0),\frac{\partial F_2}{\partial \xi}(x,0) \cdot \xi,0) 
&\mbox{ if } t=0 \\
(F_1(x,t\xi),\frac{1}{t}F_2(x,t\xi),t) &\mbox{ if } t>0
\end{array}\right.
$$
is a $\ci$ map and one proves the next proposition.
\begin{proposition}\label{gtilde}
Let $g\in \src (\Omega_{V'}^{U'})$, then 
$\tilde{g}:= g\circ \tilde{F} \in \src (\Omega_{V}^{U})$.
\end{proposition}

With the last compatibility result in hand we are ready to give the
main definition in this work.

\begin{definition}\label{src}
Let $g \in \ci (\Dnc_{X}^{M}) $.
\begin{itemize}
\item[(a)]We say that $g$ has $\it{conic\, compact\, support}$ $K$, 
if there exists a compact subset
 $K\subset M \times [0,1]$ with $K_0:=K\cap (M\times \{ 0\}) \subset X$ (conic
 compact relative to $X$) such that if $t\neq 0$ and 
$(m,t) \notin K$ then $g(m,t)=0$.
\item[(b)]We say that $g$ is $\it{rapidly\, decaying\, at\, zero}$ if for every
$(\Uo,\phi)$  $X$-slice chart
and for every $\chi \in \ci_c(\Uo \times [0,1])$, the map 
$g_{\chi}\in
\ci(\Omega_{V}^{U})$
given by
\[ g_{\chi}(x,\xi ,t)= (g\circ \varphi^{-1})(x,\xi ,t) 
\cdot (\chi \circ p \circ \varphi^{-1})(x,\xi ,t) \]
is in  $\src (\Omega_{V}^{U})$, where $p$ is the projection 
$p:\Dnc_{X}^{M} \rightarrow M\times
[0,1]$ given by $(x,\xi,0)\mapsto (x,0)$, and 
$(m,t)\mapsto (m,t)$ for $t\neq 0$, and $\varphi$ is simply the composition 
$\Dnc_{\Vo}^{\Uo}\stackrel{\tilde{\phi}}{\rightarrow}\Dnc_{V}^{U}
\stackrel{\Psi^{-1}}{\rightarrow}\Omega_{V}^{U}$. 
\end{itemize}

Finally, we denote by $\src (\Dnc_{X}^{M})$ the set of functions 
$g\in \ci(\Dnc_{X}^{M})$ that are rapidly decaying at zero 
with conic compact support.

\end{definition}

We are going to state an important property of the last
construction that we use in the next section.\\
First of all, let us recall the notion of Schwartz space associated to a
vector bundle.
\begin{definition} 
Let $(E,p,X)$ be a smooth vector bundle over a $\ci$ manifold $X$. We
define the Schwartz space $\sw (E)$ as the set of $\ci$ functions
$g \in \ci (E)$ such that $g$ is a Schwartz function at each fiber (uniformly)
and $g$ has compact support in the direction of $X$, $\it{i.e.}$, if
there exists a compact subset $K\subset X$ such that $g(E_x)=0$ for
$x\notin K$.
\end{definition}
We have the following result:

\begin{proposition}\label{sure0}
The evaluation at zero 
$\src (\Dnc_{X}^{M}) \stackrel{e_0}{\rightarrow} \sw (\Nb_{X}^{M})$ is a
surjective linear map.
\end{proposition}

\subsection{Schwartz type algebra for the Tangent groupoid}

In this section we define an algebra structure on $\src (\gr^T)$. In order to do it, 
we are going to use the functoriality of the construction 
$\src (\Dnc_{X}^{M})$ that we have seen above.\\
We start by defining a function 
$m_{r,c}: \src (\Dnc_{\go}^{\gr^{(2)}})
\rightarrow \src (\Dnc_{\gr^{(0)}}^{\gr})$ by the following formulas:\\
For $F\in \src (\Dnc_{\go}^{\gr^{(2)}})$

\[ m_{r,c}(F)(x,\xi ,0)= \int_{T_x\gr_x}F(x, \xi -\eta ,
\eta,0)
d\mu_x(\eta) \] 
and
\[ m_{r,c}(F)(\gamma,t)= \int_{\gr_{s(\gamma)}}F(\gamma \circ \delta^{-1},
\delta ,t) t^{-q}d\mu_{s(\gamma)}(\delta). \]
We have the following proposition:

\begin{proposition}\label{elproducto}
$m_{r,c}: \src ((\gr^T)^{(2)})
\rightarrow \src (\gr^T)$ is a well defined linear map.
\end{proposition}

Here we are again identifing $\Dnc_{\go}^{\gr^{(2)}}$ with
$(\gr^T)^{(2)}$ and the map above is nothing else than the integration
along the fibers of $m^T:(\gr^T)^{(2)}\rightarrow \gr^T$.\\

We are ready now to define the product in $\src (\gr^T)$.

\begin{definition}\label{conv}
Let $f,g \in \src (\gr^T)$, we define a function $f*g$ in
$\gr^T$ by 
\[ (f*g)(x,\xi,0)=\int_{T_x\gr_x}f(x, \xi -\eta ,0)
g(x,\eta,0)
d\mu_x(\eta)  \]
and
\[ (f*g)(\gamma,t)=  \int_{\gr_{s(\gamma)}}f(\gamma \circ \delta^{-1},t)
g(\delta ,t) t^{-q}d\mu_{s(\gamma)}(\delta) \]
for $t\neq 0$.
\end{definition}

\begin{proposition}
The product  $*$ is well defined and associative.
\end{proposition}

Thanks to the proposition \ref{sure0} we have an exact sequence of algebras

\begin{align}\label{se}
0 \longrightarrow J \longrightarrow \src (\gr^T)
\stackrel{e_0}{\longrightarrow} \sw (A\gr) \longrightarrow 0,
\end{align}
where $J=Ker(e_0)$ by definition.

\section{Analytic index of order $k$}

This last section is devoted to the construction of the
index announced in the introduction.\\
The main reason to construct the algebra $\src ({\gr^T})$ is that
the "Schwartz algebras" have in general the good $K-$theory
groups. For example, we are interested in the symbols of $\gr$-PDO and
more precisely in their homotopy classes in $K$-theory, that is, we are
interested in the group $K_0(A^*\gr)=K_0(C_0(A^*\gr))$. Here it
would not be enough to take the $K-$theory of $\ci_c(A\gr)$ (see for
example \cite{Co2}), however it is enough to consider the Schwartz
algebra $\sw (A^*\gr)$. Indeed, the Fourier transform shows 
that this last algebra is stable under holomorphic calculus
on $C_0(A^*\gr)$ and so it has the "good" $K$-theory, meaning that 
$K_0(A^*\gr)=K_0(\sw(A^*\gr))$.\\
By applying $K-$theory to the exact sequence (3) we obtain
\begin{align}\label{se}
K_0(J) \longrightarrow K_0(\src (\gr^T))
\stackrel{e_0}{\longrightarrow} K_0(\sw (A^*\gr)) \longrightarrow 0.
\end{align}
The surjectivity of $K_0(\src (\gr^T))
\stackrel{e_0}{\rightarrow} K_0(\sw (A^*\gr))$ comes from the existence
of a Pseudodifferential calculus on $\gr^T$.\\
Our main result is the following:

\begin{theorem}\label{teo}
There is a well defined group morphism
$$ ind_{a}^{h,k}: K^0(A^*\gr)\rightarrow K_{0}^{h,k}(\gr) $$
satisfying $ind_{a}^{h,k} \circ e_0=e_{1}^{h,k}$ where 
$e_0:K_0(\src (\gr^T)) \rightarrow K_0(\sw (A^*\gr))$ 
and $e_{1}^{h,k}$ is the composition 
$$K_0(\src (\gr^T))\stackrel{e_1}{\rightarrow} K_0(\cg) \rightarrow 
K_0(\ck) \stackrel{\pi}{\rightarrow} K_{0}^{h,k}(\gr).$$
Moreover, we have a commutative diagram

\[
\xymatrix{
\{\gr-pdo ell.\} \ar[dd]_{symb} \ar[r]^{ind}&
K_0(\ci_c(\gr)) \ar[d] & \\ & K_{0}^{h,k}(\gr) \ar[d] & \\
K_0(A^*\gr) \ar[r]_{ind_a} \ar[ru]^{ind_{a}^{h,k}}& K_0(C_{r}^{*}(\gr)) & 
.}
\]

\end{theorem}

\begin{definition}[Analytic index $ind_{a}^{h,k}$]
Let $\sigma\in K_0(A^*\gr)$ and take $w_{\sigma}\in K_0(\src (\gr^T))$ 
with $e_0(w_{\sigma})=\sigma$. The morphism given by the last theorem is explicity
calculated as  
$$ind_{a}^{h,k}(\sigma):= e_{1}^{h,k}(w_{\sigma}).$$
We call this morphism $\it{the\, analytic\, index\, of\, order}$ $k$ $\it{of}$ $\gr$. 
Putting this in a diagram, it looks as follows:

\[
\xymatrix{
K_0(J) \ar[r] &
K_0(\src (\gr^T)) \ar[d]_{e_{1}^{h,k}} \ar[r]^{e_0}
& K_0(\sw (A^*\gr)) \ar[r] \ar[ld]^{ind_{a}^{h,k}}& 0 \\
& K_{0}^{h,k}(\gr) & &
}
\]

\end{definition}

The main point of the theorem is that $ind_{a}^{h,k}$ is well defined,
because it will evidently be a group morphism. Hence, in order to prove
that it is well defined, we would like to have the following property:
\begin{itemize}
\item[(P)]If $w\in K_0(\src (\gr^T))$ with $w\in Ker(e_{0})$, then
 $e_{1}^{k}(w)\sim_h0$ in $ K_0(\ck)$.
\end{itemize}
First, consider the canonical projection $t: \gr^T \rightarrow
[0,1]$. Thanks to the condition about $\partial_t$ 
that we imposed in the definition of $\src (\gr^T)$ we have that (as
an easy consequence of the Taylor developpement)
$$J=t \cdot \src (\gr^T). $$ 
We are not going to work directly with $K_0(J)$ but with the K-theory
groups of powers $J^N$, for $N\in \Nat$.\\
Let $k\in \Nat$ and $q:=dim\, \gr_x$, we define
\[ \varphi_{k} : J^{k+q} \rightarrow \ckt \]  
by the formula (expressing $J^{k+q}=t^{k+q} \cdot \src
(\gr^T)$):\\
$$
\varphi_k (t^{k+q}\cdot f)(\gamma,t) = \left\{ 
\begin{array}{cc}
0 &\mbox{ if } t=0 \\

t^kf(\gamma,t) &\mbox{ if } t\neq 0
\end{array}\right.
$$
Then we can prove:

\begin{proposition}
With the above definition we have an algebra morphism

\[ \varphi_{k} : J^{k+q} \rightarrow \ckt \]
It satisfies $e_0 \circ \varphi_{k}=0$ by construction.
\end{proposition}

We note also by $\varphi_{k} : K_0(J^{k+q}) \rightarrow K_0(C_{c}^{k}(\gr))$
the morphism induced in  $K$-theory by $\varphi_{k}$. Now, 
using the fact that $J/J^N$ is
a nilpotent ring we check very easily that the morphism $K_0(J^N)
 \stackrel{j}{\rightarrow} K_0(J)$ induced in K-theory by the
 inclusion $J^N \hookrightarrow J$ is surjective and so we obtain the following
 commutative diagram:\\

\[
\xymatrix{
0 & & & \\
K_0(J) \ar[r] \ar[u] &
K_0(\src (\gr^T)) \ar[r]^{e_0} \ar[dr]^{e_{1}^{k}}
& K_0(\sw (A\gr)) \ar[r] & 0\\
K_0(J^{k+q}) \ar[u]^{j} \ar[r]^{\varphi_{k}} & 
K_0(C_{c}^{k}(\gr \times [0,1])) \ar[r]^{e_1} & K_0(C_{c}^{k}(\gr)) \\
& & &
}
\]

We immediately get the desired property $(P)$. 
Finally, the
commutativity of the diagram cited in the theorem is a consequence of
the possibility of doing Pseudodifferential calculus on the tangent
groupoid and the fact that the $C^*-$analytic index is also obtained
as a "deformation" through the tangent groupoid, as was shown by
Monthubert and Pierrot in \cite{MP}. 

\section{References}

\renewcommand{\refname}{}    

\vspace*{-36pt}              

\end{document}